\documentclass{amsart}
\usepackage{amsmath,amssymb,amscd,amsthm}
\newtheorem{formula}{}[section]
\newtheorem{corollary}[formula]{Corollary}
\newtheorem{lemma}[formula]{Lemma}
\newtheorem{theorem}[formula]{Theorem}
\theoremstyle{definition}
\newtheorem{definition}[formula]{Definition}

\theoremstyle{remark}

\renewcommand{\l}{\lambda}
\newcommand{\f}{\varphi}
\renewcommand{\O}{\Omega}
\newcommand{\M}{{\mathrm M}}
\newcommand{\C}{\mathbb C}
\newcommand{\R}{\mathbb R}
\newcommand{\Z}{\mathbb Z}

\newcommand{\<}{\langle}
\renewcommand{\>}{\rangle}
\newcommand{\ind}{\mathop{\rm ind}\nolimits}
\newcommand{\sign}{\mathop{\rm sign}}
\newcommand{\td}{\mathop{\rm td}}
\newcommand{\bcp}{\overline{\C P}{}^2}

\begin{document}

\title[Hirzebruch genera of manifolds with torus action]
{Hirzebruch genera of manifolds with\\ torus action}
\author{Taras E. Panov}
\thanks{This work was partly supported by the Russian Foundation for
Fundamental Research (grant no.~99-01-00090).}
\address{Department of Mathematics and Mechanics, Moscow State
University, Moscow 119899 RUSSIA}
\email{tpanov@mech.math.msu.su}

\begin{abstract}
A quasitoric manifold is a smooth $2n$-manifold $M^{2n}$ with an
action of the compact torus $T^n$ such that the action is locally
isomorphic to the standard action of $T^n$ on $\C^n$ and the orbit space is
diffeomorphic, as manifold with corners, to a simple polytope $P^n$.
The name refers to the fact that topological and
combinatorial properties of quasitoric manifolds are similar to that of
non-singular algebraic toric varieties (or toric manifolds).
Unlike toric varieties, quasitoric
manifolds may fail to be complex; however, they always admit a stably
(or weakly almost) complex structure, and their cobordism classes generate
the complex cobordism ring. As it have been recently shown by Buchstaber and
Ray, a stably complex structure on a quasitoric manifold is defined in purely
combinatorial terms, namely, by an orientation of the polytope and a function
from the set of codimension-one faces of the polytope to primitive
vectors of an integer lattice. We calculate the $\chi_y$-genus of a
quasitoric manifold with fixed stably complex structure in terms of the
corresponding combinatorial data. In particular, this gives explicit formulae
for the classical Todd genus and signature. We also relate our results with
well-known facts in the theory of toric varieties.
\end{abstract}

\maketitle

\section*{Introduction}

Manifolds with torus action arise in different areas of
topology, algebraic and differential geometry, and mathematical
physics. Specific properties of torus
action or additional structures on manifolds usually allow
to solve corresponding problems by geometrical or combinatorial methods.
The most well known examples here are
Hamiltonian torus actions in symplectic geometry and smooth toric
varieties in algebraic geometry.  Both cases allow a natural topological
generalization, namely {\it quasitoric manifolds}, introduced by Davis and
Januszkiewicz in~\cite{DJ}. (Davis and Januszkiewicz used the term ``toric
manifolds"; the term ``quasitoric manifolds" firstly appeared in~\cite{BP2}
and~\cite{BR2} because of the reasons discussed below.) A~quasitoric
manifold is a manifold with a torus action that satisfies two natural
conditions. The first one is that the action locally looks like the
standard torus action on a complex space by diagonal matrices. If this
condition is satisfied, then the orbit space is a manifold with corners;
the second condition is that this manifold with corners is diffeomorphic
to a convex simple polytope.  These two properties are well known for
smooth algebraic toric varieties~\cite{Da},~\cite{Fu}. Though quasitoric
manifolds retain most topological and combinatorial properties of smooth
toric varieties, they may fail to admit a complex structure. Like toric
varieties, quasitoric manifolds can be defined in purely combinatorial
terms. Namely, any quasitoric manifold is defined by combinatorial data:
the lattice of faces of a simple polytope and a {\it characteristic
function} that assigns an integer primitive vector defined up to sign to
each facet.  Despite their simple and specific definition, quasitoric
manifolds in many cases may serve as model examples (for instance, as it
was shown in~\cite{BR1}, each complex cobordism class contains a
quasitoric manifold). All these facts allow to use quasitoric manifolds
for solving topological problems by combinatorial methods and vice versa.
A number of such relations was firstly discovered in the theory of toric
varieties.  Some applications were obtained in~\cite{BP1}, \cite{BP2},
where quasitoric manifolds are studied in the general context of
``manifolds defined by simple polytopes". Another example of
interplay between topology and combinatorics is the calculation of the
$KO$-theory of quasitoric manifolds held in~\cite{BB}.

Besides, it is important to mention that
the term ``quasitoric manifold" is not common. Some authors starting from
Davis and Januszkiewicz call the object of our study just ``toric manifolds".
However, we prefer to call those manifolds ``quasitoric", reserving the term
``toric manifold" for smooth toric varieties, as commonly used in
algebraic geometry (see, for instance,~\cite{Ba}).

In the present paper we calculate
well-known cobordism invariant, the $\chi_y$-genus, for quasitoric manifolds
in terms of the corresponding combinatorial data. The most important
particular cases here are the signature and the Todd genus, which correspond
to the values $y=1$ and $y=0$ respectively. The signature is an oriented
cobordism invariant and is defined for any oriented manifold (it equals zero
in dimensions other than $4k$). At the same time the Todd genus, as well
as the general $\chi_y$-genus, requires the Chern classes of manifold to be
defined. As it is mentioned above, a quasitoric manifold is not
necessarily complex, however as it was recently shown by Buchstaber
and Ray~\cite{BR2}, it always admits a stably (or weakly almost) complex
structure. Moreover, stably complex structures (i.e. complex structures in
the stable tangent bundle) on quasitoric manifolds are also defined
combinatorially, namely, by specifying an orientation of the simple
polytope and choosing signs for vectors corresponding to facets. This in
turn is equivalent to a choice of orientations for the manifold and
all submanifolds corresponding to facets. A quasitoric manifold with such
additional structure was called in~\cite{BR2} {\it multioriented}. Hence,
a multioriented quasitoric manifold determines a complex cobordism class,
for which one can define characteristic numbers and complex Hirzebruch
genera. We calculate the $\chi_y$-genus of a multioriented (i.e. with
fixed stably complex structure) quasitoric manifold in terms of its
combinatorial data. To do this, we construct a circle action with only
isolated fixed points and then apply Atiyah and Bott's generalized
Lefschetz fixed point theorem~\cite{AB}. The obtained formula allows to
calculate the $\chi_y$-genus as a sum of contributions corresponding to
the vertices of polytope. These contributions depend only on the ``local
combinatorics" near the vertex. Our formula contains one external
parameter (an integer primitive vector $\nu$), which does not affect the
answer. In the particular case of smooth toric varieties the Todd genus
(or the arithmetic genus) is always equal to 1~\cite{Fu}. However, this
fact fails to be true for general quasitoric manifolds, since the Todd
genus is a complex cobordism invariant and each cobordism class contains a
quasitoric representative~\cite{BR1}.

\section{Quasitoric manifolds and
polytopes, facet and edge vectors, and stably complex structures}

Here we briefly discuss the definition of quasitoric manifolds and introduce
stably complex structures on them, following~\cite{DJ} and~\cite{BR2}. We
also describe some new relations between the combinatorial data associated
to a quasitoric manifold.

A {\it convex polytope} $P^n$ of dimension $n$ is a bounded set in
$\R^n$ that is obtained as the intersection of a finite number of
half-spaces:
\begin{equation}
\label{ptope}
  P^n=\{x\in\R^n:\<l_i,x\>\ge-a_i,\; i=1,\ldots,m\}
\end{equation}
for some $l_i\in(\R^n)^{*}$, $a_i\in\R$.
A convex polytope $P^n$ is called {\it simple} if its bounding
hyperplanes are in general position at each vertex, i.e. there exactly $n$
codimension-one faces (or {\it facets}) meet at each vertex. It follows that
any point of a simple polytope
lies in a neighbourhood that is affinely isomorphic to
an open set of the positive cone
$$
  \R^n_+=\{(x_1,\ldots,x_n)\in\R^n:x_1\ge0,\ldots,x_n\ge0\}.
$$
This means exactly that a simple polytope $P^n$ is an $n$-dimensional
{\it manifold with corners}. The faces of $P^n$ of all dimensions
form a partially ordered set with respect to inclusion, which is
called the {\it lattice of faces} of $P^n$. We say that two polytopes
are {\it combinatorially equivalent} if they have same lattices of faces.
Two polytopes are combinatorially equivalent if and only if they are
diffeomorphic as manifolds with corners.

Let $M^{2n}$ be a compact $2n$-dimensional manifold with an action of the
compact torus $T^n$. One can view $T^n$ as a subgroup of the complex torus
$(\C^{*})^n$ in the standard way:
$$
  T^n=\{\bigl(e^{2\pi i\f_1},\ldots,e^{2\pi i\f_n}\bigr)\in\C^n\},
$$
where $(\f_1,\ldots,\f_n)$ runs over $\R^n$.
We say that a $T^n$-action is {\it
locally isomorphic to the standard} diagonal action of $T^n$ on $\C^n$ if
every point $x\in M^{2n}$ lies in some $T^n$-invariant neighbourhood $U(x)$
for which there exists an equivariant homeomorphism $f:U(x)\to W$ with some
($T^n$-stable) open subset $W\subset\C^n$. The last statement means that
there is an automorphism $\theta:T^n\to T^n$ such that $f(t\cdot
y)=\theta(t)f(y)$ for all $t\in T^n$, $y\in U(x)$. The orbit space for such
an action of $T^n$ is an $n$-dimensional manifold with corners; we refer to
$M^{2n}$ as a {\it quasitoric manifold} if the orbit space is
diffeomorphic, as manifold with corners, to a simple polytope $P^n$.
Thus, the orbit space of a quasitoric manifold is decomposed into faces in
such a way that points from the relative interior of each $k$-face
correspond to orbits with same isotropy subgroup of codimension $k$. In
particular, the action of $T^n$ is free over the interior of $P^n$, while
the vertices of $P^n$ correspond to $T^n$-fixed points of $M^{2n}$.

Now let $M^{2n}$ be a quasitoric manifold with orbit space $P^n$,
and let $\mathcal F=\{F_1,\ldots, F_m\}$ be the set of
codimension-one faces (facets) of the polytope $P^n$, $m=\sharp\mathcal F$.
The interior of each facet $F_i$ consists of orbits with the same
one-dimensional isotropy subgroup $G_{F_i}$. This one-parameter
subgroup of $T^n$ is determined by an integer primitive vector
$\l_i=(\l_{1i},\ldots,\l_{ni})^\top$ in the
corresponding lattice $L\simeq\Z^n$:
\begin{equation}
\label{fisotr}
  G_{F_i}=\{\bigl(e^{2\pi i\l_{1i}\f},\ldots,e^{2\pi
  i\l_{ni}\f}\bigr)\in T^n\},
\end{equation}
where $\f\in\R$. Of course, the vector $\l_i$ is defined only up to
sign. In this way one defines the {\it characteristic function}
$\l:\mathcal F\to\Z^n$ that takes a facet to the corresponding primitive
vector. By the definition of quasitoric manifolds, the characteristic
function satisfies the following condition: if $n$ facets
$F_{i_1},\ldots,F_{i_n}$ meet at same vertex $p$, i.e.
$p=F_{i_1}\cap\cdots\cap F_{i_n}$, then the integer vectors
$\l(F_{i_1}),\ldots,\l(F_{i_n})$ constitute an integer basis of the
lattice $L\simeq\Z^n$. As it was shown in~\cite{DJ}, the characteristic
pair $\left(P^n,\l:\mathcal F\to\Z^n\right)$ satisfying the above
condition defines a quasitoric manifold $M^{2n}$ uniquely up to an
equivariant diffeomorphism. Once the numeration of facets of $P^n$ is
fixed, the characteristic function can be viewed as the integer $(n\times
m)$-matrix $\Lambda$ whose $i$-th column is the vector $\l(F_i)$. Each
vertex $p$ of $P^n$ can be represented as the intersection of $n$ facets:
$p=F_{i_1}\cap\cdots\cap F_{i_n}$; we denote by
$\Lambda_{(p)}=\Lambda_{(i_1,\ldots,i_n)}$ the minor matrix of $\Lambda$
formed by the columns $i_1,\ldots,i_n$. It follows from the above that
\begin{equation}
\label{L}
  \det\Lambda_{(p)}=\det(\l_{i_1},\ldots,\l_{i_n})=\pm1.
\end{equation}
In the sequel we refer to the vectors $\l_i=\l(F_{i})$ as {\it facet
vectors}. Note again that for now they are defined only up to sign.

The next things we need in order to define genera of a quasitoric manifold
is an orientation and stably complex structures. We fix an orientation of
$T^n$ and specify an orientation of the polytope $P^n$ by orienting the
ambient space $\R^n$. These orientation data define an orientation of the
quasitoric manifold $M^{2n}$. The above construction of the characteristic
function $\l$ involves arbitrary choice of signs for the vectors
$\l(F_i)$. The inverse image $\pi^{-1}(F_i)$ of the facet $F_i$ under the
orbit map $\pi:M^{2n}\to P^n$ is a quasitoric submanifold $M_i\subset
M^{2n}$ of codimension 2. This facial submanifold $M_i$ belongs to the
fixed point set of the circle subgroup $G_{F_i}\subset T^n$
(see~(\ref{fisotr})), which therefore acts on the normal bundle
$\nu_i:=\nu(M_i\subset M^{2n})$.  Hence, one can invest $\nu_i$ with a
complex structure (and orientation) by specifying a sign of the
primitive vector $\l(F_i)$. Given an orientation of $M^{2n}$ (i.e. an
orientation of $P^n$) and an orientation of $\nu_i$ (i.e. a sign of
$\lambda_i$), one can orient the facial submanifold $M_i$ (or,
equivalently, the facet $F_i\subset P^n$). Now, the oriented submanifold
$M_i\subset M^{2n}$ of codimension 2 gives rise, by the standard
topological construction, to a complex line bundle $\sigma_i$ over
$M^{2n}$ that restricts to $\nu_i$ over $M_i$. The following theorem was
proved by Buchstaber and Ray in~\cite{BR2}:

\begin{theorem}
\label{stcomst}
  The following isomorphism of real oriented $2m$-bundles holds for any
  quasitoric manifold $M^{2n}$:
  $$
    \tau(M^{2n})\oplus\R^{2(m-n)}\simeq\sigma_1\oplus\cdots\oplus\sigma_m.
  $$
  Here $\tau(M^{2n})$ is the tangent bundle and $\R^{2(m-n)}$ denotes the
  trivial $2(m-n)$-bundle over $M^{2n}$.
\end{theorem}

\noindent An oriented quasitoric manifold $M^{2n}$ together with fixed
orientations of facial submanifolds $M_i=\pi^{-1}(F_i)$ was called
in~\cite{BR2} {\it multioriented}. It follows that signs for the facet
vectors $\{\l_i\}$ of a multioriented quasitoric manifold are defined
unambiguously. Theorem~\ref{stcomst} shows that a multioriented
quasitoric manifold can be invested with a canonical stably complex structure.
Thus, an oriented simple polytope $P^n$ with characteristic matrix
$\Lambda$ not only define a (multioriented) quasitoric manifold,
but also specify a cobordism class in the complex cobordism ring $\O_U$.

Given an edge (1-dimensional
face) $E_j$ of $P^n$, we see that points from its relative interior in $P^n$
correspond to orbits with the same $(n-1)$-dimensional isotropy subgroup
$G_{E_j}$. This subgroup can be written as
\begin{equation}
\label{eisotr}
  G_{E_j}=\{\bigl(e^{2\pi i\f_1},\ldots,e^{2\pi i\f_n}\bigr)\in T^n:
  \mu_{1j}\f_1+\ldots+\mu_{nj}\f_n=0\},
\end{equation}
i.e. it is determined by a primitive (co)vector
$\mu_j=(\mu_{1j},\ldots,\mu_{nj})^\top$ in the dual (or weight) lattice
$W=L^{*}$. We refer to this $\mu_j$ as {\it edge vector}; again it is
defined only up to sign. The edge vectors satisfy the condition similar to
that for the facet vectors: if $E_{j_1},\ldots,E_{j_n}$ are edges that meet at
same vertex $p$, then
\begin{equation}
\label{M}
  \det\M_{(p)}=\det(\mu_{j_1},\ldots,\mu_{j_n})=\pm1.
\end{equation}
Here $\M_{(p)}$ is the square matrix with columns
$\mu_{j_1},\ldots,\mu_{j_n}$. The above equality means that the vectors
$\mu_{j_1},\ldots,\mu_{j_n}$ constitute a basis of $W\simeq\Z^n$.

The following lemma enables to choose signs of edge vectors for a
multioriented quasitoric manifold unambiguously ``locally" at each vertex.
\begin{lemma}
\label{lm}
  Given a vertex $p$ of $P^n$, one can choose signs for edge vectors
  $\mu_{j_1},\ldots,\mu_{j_n}$ meeting at $p$ in
  such a way that
  $$
   \M_{(p)}^\top\cdot\Lambda_{(p)}=E,
  $$
  where $E$ denotes the identity matrix, and
  the matrices $\Lambda_{(p)}$, $\M_{(p)}$ are those
  from~{\rm (\ref{L})},~{\rm (\ref{M})}.
\end{lemma}
\begin{proof}
Since we are interested only in facet and edge vectors meeting at the
vertex $p=F_{i_1}\cap\cdots\cap F_{i_n}$, we may renumerate the
edge vectors $\mu_{j_1},\ldots,\mu_{j_n}$ at $p$ by the index set
$\{i_1,\ldots,i_n\}$
of the facet vectors at $p$. To do this we just set $j_k=i_k$ if
$E_{j_k}\not\subset F_{i_k}$ (i.e. a facet vector and an edge vector have
same index if the corresponding facet and edge of $P^n$ span the whole
$\R^n$). Then for $k\ne l$ one has $E_{i_k}\subset F_{i_l}$. Hence,
$G_{F_{i_l}}\subset G_{E_{i_k}}$ and
\begin{equation}
\label{mlne}
  \<\mu_{i_k},\l_{i_l}\>=0, \quad k\ne l,
\end{equation}
(see~(\ref{fisotr}) and~(\ref{eisotr})). Now, since $\mu_{i_k}$ is a primitive
vector, it follows from~(\ref{mlne}) that $\<\mu_{i_k},\l_{i_k}\>=\pm1$.
Changing the sign of $\mu_{i_k}$ if necessary, we obtain
\begin{equation}
\label{mle}
  \<\mu_{i_k},\l_{i_k}\>=1,
\end{equation}
which together with~(\ref{mlne}) gives $\M_{(p)}^\top\cdot\Lambda_{(p)}=E$,
as needed.
\end{proof}

Throughout the rest of this paper we assume that once a stably complex
structure of quasitoric manifold is fixed (i.e. signs for
columns of the characteristic $(n\times m)$-matrix $\Lambda$ are specified),
signs of edge vectors for each vertex are chosen as in
Lemma~\ref{lm}.

We fix an orientation of the torus $T^n$ once and forever; then
a choice of orientation for $M^{2n}$ is equivalent to a choice of
orientation for the polytope $P^n\subset\R^n$. Hence, edges
of $P^n$ meeting at the same vertex $p$ can be ordered canonically in such a
way that the ordered set of vectors along the edges pointing out of $p$
constitute a positively oriented basis of $\R^n$ (i.e. an orientation of
$P^n$ defines an ordering of edges at each vertex). In the sequel we assume
that once an orientation of $M^{2n}$ is fixed, the edge vectors
$\mu_{i_1},\ldots,\mu_{i_n}$ meeting at $p$ are ordered
in accordance with the above ordering of edges at $p$.
In this situation the edge vectors $\mu_{i_1},\ldots,\mu_{i_n}$ themselves
may constitute a positively or negatively oriented basis of $\R^n$. (We
assume that signs for $\mu_{i_1},\ldots,\mu_{i_n}$ are determined by
Lemma~\ref{lm}.) Thus, we come to the following
\begin{definition}
\label{sign}
Given a multioriented quasitoric manifold
$M^{2n}$ with orbit space $P^n$, the {\it sign} of a vertex $p\in P^n$ is
$$
  \sigma(p)=\det\M_{(p)}=\det(\mu_{i_1},\ldots,\mu_{i_n}),
$$
where $\mu_{i_1},\ldots,\mu_{i_n}$ are canonically ordered
edge vectors meeting at $p$.
\end{definition}
Obviously, Lemma~\ref{lm} shows that
$$
  \sigma(p)=\det\Lambda_{(p)}=\det(\l_{i_1},\ldots,\l_{i_n}),
$$
where $\l_{i_1},\ldots,\l_{i_n}$ are the facet vectors at $p$.
We mention that the
above definition of the sign of vertex already appeared in~\cite{Do}
while studying characteristic functions and quasitoric manifolds over
given simple polytope.

\medskip

As it is mentioned in the introduction, quasitoric manifolds can be viewed
as a topological generalization of algebraic smooth {\it toric varieties}.
Toric variety~\cite{Da} is a normal algebraic variety on which the complex
torus $(\C^{*})^n$ acts with a dense orbit. Though any smooth toric variety
(also called {\it toric manifold}) is
quasitoric manifold, we illustrate our above constructions in the more
restricted case of smooth projective toric varieties arising from
polytopes. Any such toric variety is obtained via the standard procedure
(see~\cite{Fu}) from a simple polytope~(\ref{ptope}) with vertices
in the integer lattice $\Z^n\subset\R^n$ (we refer to such polytope as
{\it integral}). Normal covectors $l_i$ of facets of an integral
polytope $P^n$ can be chosen integer and primitive. The toric variety $M_P$
obtained from an integral simple polytope $P^n$ is necessarily projective and
has complex dimension $n$. The compact torus $T^n\subset(\C^{*})^n$ acts on
$M_P$ with orbit space $P^n$; moreover, there is a map from $M_P$ to $\R^n$
which is constant on $T^n$-orbits and has image $P^n$ (the {\it moment
map}).  The variety $M_P$ is smooth whenever for each vertex $p\in P^n$ the
normal covectors $l_{i_k}$, $k=1,\ldots,n$, of facets containing $p$
constitute an integer basis of the dual lattice. This means exactly that
the map $F_i\to l_i$ defines a characteristic function. It can be easily
seen that it is exactly the characteristic function of $M_P$ regarded as a
quasitoric manifold (or, conversely, the quasitoric manifold corresponding
to the above characteristic function is equivariantly diffeomorphic to
$M_P$).  Hence, facet vectors for $M_P$ are just normal (co)vectors for
the corresponding polytope $P^n$. Edge vectors are primitive integer
vectors along the edges of $P^n$. In short, this means that in the case of
toric $M_P$ ``facet vectors" are ``vectors normal to facets" and ``edge
vectors" are ``vectors along edges". Lemma~\ref{lm} in this case says that
if an edge is contained in a facet, then the vector along the edge is
orthogonal to the normal vector of the facet, while if an edge is not
contained in a facet, then signs of the corresponding vectors can be
chosen in such a way that their scalar product is equal to 1.  It can be
shown that the canonical complex structure on a toric variety $M_P$,
regarded as a stably complex structure on a quasitoric manifold,
corresponds to orienting the facet vectors $l_i$ in such a way that all of
them are ``pointing inside the polytope $P^n$". As it follows from
Lemma~\ref{lm}, edge vectors locally should be oriented in such a way that
all of them are ``pointing out of the vertex". Note that globally
Lemma~\ref{lm} provides two signs for an edge, one for each of its ends.
These signs are always different in the case of toric varieties, however this
is not true in general. Note also that since an edge vector for a toric
variety is a vector along the edge pointing out of the vertex, one has
$\sigma(p)=1$ (see Definition~\ref{sign}) for any vertex $p$.

\section{Circle action with isolated fixed points on a quasitoric
manifold}

In this section we show that there exists a circle subgroup of $T^n$ that
acts on a quasitoric manifold with only isolated fixed points
corresponding to vertices of the underlying simple polytope. This
action will be used in the next section for calculating the $\chi_y$-genus
via contributions of fixed points.

So, we start with a multioriented quasitoric manifold $M^{2n}$ with orbit
space $P^n$ and edge vectors $\{\mu_i\}$. Hence,
$M^{2n}$ is endowed with a stably complex structure, as described in
the previous section.

\begin{theorem}
\label{sa}
Suppose that $\nu\in\Z^n$ is an integer primitive vector such that
$\<\mu_i,\nu\>\ne0$ for all edge vectors $\mu_i$. Then the circle subgroup
$S^1\subset T^n$ defined by $\nu$ acts on $M^{2n}$ with isolated fixed points
corresponding to vertices of $P^n$. In the tangent space $T_pM^{2n}$ at fixed
point corresponding to the vertex $p=F_{i_1}\cap\cdots\cap F_{i_n}$
this action induces a representation of $S^1$ with weights
$\<\mu_{i_1},\nu\>,\ldots,\<\mu_{i_n},\nu\>$.
\end{theorem}

\begin{proof}
Let us consider the action of $T^n$ near fixed point $p\in M^{2n}$
corresponding to a vertex
$p=F_{i_1}\cap\cdots\cap F_{i_n}$ of $P^n$.
Let $\mu_{i_1},\ldots,\mu_{i_n}$ be the edge vectors at $p$.
The $T^n$-action induces a unitary $T^n$-representation in the tangent space
$T_pM^{2n}$. We choose complex
coordinates $(x_1,\ldots,x_n)$ in $T_pM^{2n}$ such
that the tangent space of two-dimensional submanifold
$\pi^{-1}(E_{i_k})\subset M^{2n}$ at $p$ is given by equations
$x_1=\ldots=\widehat{x_k}=\ldots=x_n=0$ ($x_k$ is dropped).
Then, the corresponding isotropy subgroup $G_{E_{i_k}}$ is
given by the equation $\<\mu_{i_k},\f\>=0$ in $T^n$ (see~(\ref{eisotr})).
Hence, the weights of the $T^n$-representation are
$\mu_{i_1},\ldots,\mu_{i_n}$, i.e.
an element
$t=\bigl(e^{2\pi i\f_1},\ldots,e^{2\pi i\f_n}\bigr)\in T^n$
acts on $T_pM^{2n}$ as
\begin{align}
\label{taction}
  t\cdot(x_1,\ldots,x_n)&=\bigl(
  e^{2\pi i(\mu_{1i_1}\f_1+\cdots+\mu_{ni_1}\f_n)}x_1,\ldots,
  e^{2\pi i(\mu_{1i_n}\f_1+\cdots+\mu_{ni_n}\f_n)}x_n \bigr)\\
  &=\bigl(e^{2\pi i\<\mu_{i_1},\Phi\>}x_1,\ldots,
  e^{2\pi i\<\mu_{i_n},\Phi\>}x_n \bigr),\notag
\end{align}
where $(x_1,\ldots,x_n)\in T_pM^{2n}$, $\Phi=(\f_1,\ldots,\f_n)$.

A primitive vector $\nu=(\nu_1,\ldots,\nu_n)^\top\in\Z^n$ defines a
one-parameter
circle subgroup $\{\bigl(e^{2\pi i\nu_1\f},\ldots,e^{2\pi i\nu_n\f}\bigr),
\f\in\R\}\subset T^n$. It follows from~(\ref{taction}) that this circle acts
on $T_pM^{2n}$ as
$$
  s\cdot(x_1,\ldots,x_n)=
  \bigl(e^{2\pi i\<\mu_{i_1},\nu\>\f}x_1,\ldots,
  e^{2\pi i\<\mu_{i_n},\nu\>\f}x_n \bigr),
$$
where $s=e^{2\pi i\f}\in S^1$.
The fixed point $p$ is isolated if all the weights
$\<\mu_{i_1},\nu\>,\ldots,$ $\<\mu_{i_n},\nu\>$ of the $S^1$-action are
non-zero.  Thus, if $\<\mu_i,\nu\>\ne0$ for all edge vectors, then the
$S^1$-action on $M^{2n}$ defined by $\nu$ has only isolated fixed points.
\end{proof}

Note that the condition $\<\mu_i,\nu\>\ne0$ for all $\mu_i$ from the above
theorem means that the primitive vector $\nu$ is ``of general position". In
the case of smooth toric variety constructed from integral simple polytope
$P^n$ the above condition is that the vector $\nu$ is not orthogonal to any
edge of $P^n$ (or, equivalently, no hyperplane normal to $\nu$
intersects $P^n$ by a face of dimension $>0$).

In the next section we will need the following
\begin{definition}
\label{ind}
Given the $S^1$-action on $M^{2n}$ defined by a primitive vector $\nu$, the
{\it index} of the vertex $p=F_{i_1}\cap\cdots\cap F_{i_n}$ is
$$
  \ind_\nu(p)=\{\sharp k:\<\mu_{i_k},\nu\><0\},
$$
i.e. $\ind_\nu(p)$ equals the number of negative weights at $p$.
\end{definition}

The following lemma shows that the index of any vertex $p$ can be defined
directly by means of facet vectors at $p$ without calculating edge vectors.

\begin{lemma}
  Let $p=F_{i_1}\cap\cdots\cap F_{i_n}$ be a vertex of $P^n$.
  Then the index $\ind_\nu(p)$ equals the number of negative coefficients in
  the representation of $\nu$ as a linear combination of basis vectors
  $\l_{i_1},\ldots,\l_{i_n}$.
\end{lemma}
\begin{proof}
It follows from Lemma~\ref{lm} that for any vertex
$p=F_{i_1}\cap\cdots\cap F_{i_n}$ one can write
$$
  \nu=\<\mu_{i_1},\nu\>\l_{i_1}+\ldots+\<\mu_{i_n},\nu\>\l_{i_n}.
$$
Then our lemma follows from the definition of $\ind_\nu(p)$.
\end{proof}

\section{Formulae for $\chi_y$-genus, signature and Todd genus}

Here we calculate the $\chi_y$-genus of a multioriented quasitoric manifold
in terms of its characteristic pair $(P^n,\Lambda)$.  This is done by
applying the Atiyah--Hirzebruch formula for the $S^1$-action constructed in
the previous section.  We also stress upon the most important particular
cases of signature and Todd genus.

\begin{theorem}
\label{chi}
  Let $M^{2n}$ be a multioriented quasitoric manifold, and let
  $\nu\in\Z^n$ be an integer primitive vector such that
  $\<\mu_i,\nu\>\ne0$ for all edge vectors $\mu_i$. Then
  $$
    \chi_y(M^{2n})=\sum_{p\in P^n}(-y)^{\ind_\nu(p)}\sigma(p),
  $$
  where the sum is taken over all vertices of $P^n$, $\sigma(p)$ and
  $\ind_\nu(p)$ are as in definitions~{\rm \ref{sign}}
  and~{\rm \ref{ind}.}
\end{theorem}
\begin{proof}
The Atiyah--Hirzebruch formula (\cite{AH}; see also~\cite{Kr}, where it was
deduced within the cobordism theory) states that the $\chi_y$-genus of a
stably complex manifold $X$ with a $S^1$-action can be calculated as
\begin{equation}
\label{AH}
  \chi_y(X)=\sum_i(-y)^{n(F_i)}\chi_y(F_i),
\end{equation}
where the sum is taken over all $S^1$-fixed submanifolds $F_i\subset X$,
and $n(F_i)$ denotes the number of negative weights of the
$S^1$-representation in the normal bundle of $F_i\subset X$. In our case
all fixed submanifolds are isolated fixed points
corresponding to vertices $p\in P^n$. Hence, $\chi_y(F_i)=\chi_y(p)=\pm1$
depending on whether or not the orientation of $\R^n$ defined by
$\mu_{i_1},\ldots,\mu_{i_n}$ coincides with that defined by $P^n$. Thus,
for quasitoric $M^{2n}$ we may substitute $\sigma(p)$ for $\chi_y(F_i)$
in~(\ref{AH}). Theorem~\ref{sa} shows that the weights of induced
$S^1$-representation in $T_pM^{2n}$ equal
$\<\mu_{i_1},\nu\>,\ldots,\<\mu_{i_n},\nu\>$, therefore $n(F_i)$
in~(\ref{AH}) is exactly $\ind_\nu(p)$ (see Definition~\ref{ind}), and the
required formula follows.
\end{proof}

The $\chi_y$-genus $\chi_y(M^{2n})$ at $y=-1$ equals the $n$th Chern number
$c_n[M^{2n}]$ for any $2n$-dimensional stably complex manifold $M^{2n}$.
Theorem~\ref{chi} gives
\begin{equation}
\label{cn}
  c_n[M^{2n}]=\sum_{p\in P^n}\sigma(p).
\end{equation}
If $M^{2n}$ is a complex manifold (e.g., $M^{2n}$ is a smooth
toric variety) one has $\sigma(p)=1$ for all vertices $p\in P^n$ and
$c_n[M^{2n}]$ equals the Euler number $e(M^{2n})$. Hence, for
complex $M^{2n}$ the Euler number equals the number of vertices of
$P^n$ (which is well known for toric varieties). For general quasitoric
$M^{2n}$ the Euler number is also equal to the number of vertices of
$P^n$ (since the Euler number of any $S^1$-manifold equals the
sum of Euler numbers of fixed submanifolds), however this number may
be different from $c_n[M^{2n}]$.

\medskip

The $\chi_y$-genus at $y=1$ equals the signature (or the
$L$-genus). Theorem~\ref{chi} gives in this case
\begin{corollary}
\label{signat}
 The signature of a multioriented quasitoric manifold $M^{2n}$
 can be calculated as
 $$
   \sign(M^{2n})=\sum_{p\in P^n}(-1)^{\ind_\nu(p)}\sigma(p),
 $$
 where the sum is taken over all vertices of $P^n$.
\end{corollary}

Being an invariant of an oriented cobordism class the signature of a
quasitoric manifold $M^{2n}$ does not depend on a stably complex structure
(i.e. on a choice of signs for facet vectors)
and is determined only by an orientation of $M^{2n}$ (or~$P^n$).
This can be seen directly from our above considerations. Indeed, we have
chosen signs of edge vectors locally at each vertex as described in
Lemma~\ref{lm}.
Then Corollary~\ref{signat} states that the signature can be calculated as
the sum over vertices of values $(-1)^{\ind_\nu(p)}\sigma(p)$, where
$\ind_\nu(p)$ is the number of negative scalar products $\<\mu_{i_k},\nu\>$
and $\sigma(p)=\det(\mu_{i_1},\ldots,\mu_{i_n})$ determines the difference
between orientations of
$T_pM^{2n}$ defined by the $T^n$-representation with weights
$\mu_{i_1},\ldots,\mu_{i_n}$ and the standard $T^n$-representation.
Instead of this, we may choose signs for edge vectors in such a way
that {\it all\/} scalar products with $\nu$ are positive.
Then one obviously has
$$
  (-1)^{\ind_\nu(p)}\sigma(p)=
  \det(\tilde\mu_{i_1},\ldots,\tilde\mu_{i_n}),
$$
where $\tilde\mu_{i_k}=\pm\mu_{i_k}$ and $\<\tilde\mu_{i_k},\nu\>>0$,
$k=1,\ldots,n$. Now, the right hand side of the above equality is
independent of particular choice of signs for facet vectors (i.e.
independent of a stably complex structure). Thus, we come to the following
\begin{corollary}
\label{signor}
  The signature of an oriented quasitoric manifold $M^{2n}$ can be calculated
  as
  $$
  \sign(M^{2n})=\sum_{p\in P^n}\det(\tilde\mu_{i_1},\ldots,\tilde\mu_{i_n}),
  $$
  where $\{\tilde\mu_{i_k}\}$ are edge vectors at $p$
  with signs chosen to satisfy $\<\tilde\mu_{i_k},\nu\>>0$, $k=1,\ldots,n$.
\end{corollary}

In the case of smooth toric variety $M_P$ we have $\sigma(p)=1$ for all
vertices $p\in P^n$, and Corollary~\ref{signat} gives
\begin{equation}
\label{toricsign}
  \sign(M_P)=\sum_{p\in P^n}(-1)^{\ind_\nu(p)}.
\end{equation}
The scalar product with vector $\nu$ can be viewed as an analog of
Morse height function on the polytope $P^n\subset\R^n$
(see~\cite{Br}, \cite{Kh}, and~\cite{DJ}).  Since $\nu$ is
not orthogonal to any edge of $P^n$, we may orient the edges of $P^n$ in
such a way that the scalar product with $\nu$ increases along the edges.
Then the index of vertex $p$ is just the number of edges pointing in
(i.e. towards $p$). An easy combinatorial argument~\cite[p.~115]{Br} shows
that the number of vertices with exactly $k$ edges pointing in equals
$h_k$. Here $(h_0,h_1,\ldots,h_n)$ is the so-called {\it $h$-vector} of
the polytope $P^n$. It is defined from the equation
$$
  h_0t^n+\ldots+h_{n-1}t+h_n=(t-1)^n+f_0(t-1)^{n-1}+\ldots+f_{n-1},
$$
where $f_k$ is the number of faces of $P^n$ of codimension $(k+1)$. It is well
known~\cite{Fu} that $h_k$ equals $2k$-th Betti number of the toric variety
$M_P$ (this is also true for a general quasitoric manifold~\cite{DJ},
however we do not use this fact here). Now, using
formula~(\ref{toricsign}) we obtain that the signature of a smooth toric
variety $M_P$ is
$$
  \sign(M_P)=\sum_{k=1}^n(-1)^kh_k.
$$
This formula can be also deduced from the Hodge structure in the
cohomology of $M_P$. However, we believe that obtaining this result from our
general formula for quasitoric manifolds could be of interest.

\medskip

The next important particular case of the $\chi_y$-genus is the Todd genus
corresponding to the value $y=0$. In this case summands in the formula
from Theorem~\ref{chi} are not defined for those vertices having index 0 ,
so it requires some additional analysis.

\begin{theorem}
\label{todd}
  The Todd genus of a multioriented quasitoric manifold
  can be calculated as
  $$
    \td(M^{2n})=\mathop{\sum_{p\in P^n}}\limits_{\ind_\nu(p)=0}\sigma(p)
  $$
  (the sum is taken over all vertices of index 0). Here $\nu$ is any
  primitive vector such that $\<\mu_i,\nu\>\ne0$ for all $\mu_i$.
\end{theorem}
\begin{proof}
For each vertex $p=F_{i_1}\cap\cdots\cap F_{i_n}$ of $P^n$ the
stably complex structure on $M^{2n}$ determined by $(P^n,\Lambda)$ defines
a complex structure on $T_pM^{2n}$ via the isomorphism
\begin{equation}
\label{ts}
  T_pM^{2n}\simeq\sigma_{i_1}\oplus\cdots\oplus\sigma_{i_n}
\end{equation}
(see Theorem~\ref{stcomst}). The $S^1$-action on $M^{2n}$ determined by
$\nu\in\Z^n$ induces the complex $S^1$-representation on $T_pM^{2n}$ with
weights $w_1(p)=\<\mu_{i_1},\nu\>,\ldots,w_n(p)=\<\mu_{i_n},\nu\>$ (signs
of edge vectors are determined by Lemma~\ref{lm}). Atiyah and Bott's
generalized Lefschetz fixed point formula
(\cite{AB}, see also~\cite{HBY}) gives the following
expression for the {\it equivariant $\chi_y$-genus} of $M^{2n}$
\begin{equation}
\label{chieq}
  \chi_y(q,M^{2n})=\sum_{p\in P^n}\prod_{i=1}^k
  \frac{1+yq^{w_k(p)}}{1-q^{w_k(p)}}\sigma(p),
\end{equation}
where $q=e^{2\pi i\varphi}\in S^1$, and
$\sigma(p)=\det(\mu_{i_1},\ldots,\mu_{i_n})=\pm1$ depending on whether or
not the orientation of $T_pM^{2n}$ defined by~(\ref{ts}) coincides with that
defined by the original orientation of $M^{2n}$. (We note again that this sign
$\sigma(p)$ is equal 1 for all vertices if $M^{2n}$ is a true complex
manifold, e.g., a smooth toric variety.) Atiyah and
Hirzebruch's theorem~\cite{AH}
states that the above expression for $\chi_y(q,M^{2n})$
is independent of $q$ and equals $\chi_y(M^{2n})$. Taking the limit of the
right hand side of~(\ref{chieq}) as $q\to0$, one obtains the
Atiyah--Hirzebruch formula~(\ref{AH}) (since the $\lim_{q\to0}$ of each
summand in~(\ref{chieq}) is exactly $(-y)^{\ind(p)}\sigma(p)$). In the case
$y=0$ corresponding to the Todd genus the same limit for the summand
corresponding to a vertex $p$ equals 0 if there is at least one $w_k(p)<0$
and equals 1 otherwise. This is exactly what is stated in the theorem.
\end{proof}

As for the Todd genus of a smooth toric variety, it is easy to see that
there is only one vertex of index 0 in this case. Indeed, if we orient
edges of the polytope by means of the scalar product with $\nu$ (see the
above considerations concerning the signature of toric varieties), then
only one ``bottom" vertex will have all edges pointing out, i.e. index 0.
Since one has $\sigma(p)=1$ for all vertices $p\in P^n$,
Theorem~\ref{todd} gives $\td(M_P)=1$, which is well known (see
e.g.~\cite{Fu}).

\medskip

At the end we consider some examples illustrating our results. All our
examples are multioriented four-dimensional quasitoric manifolds.
Thus, all facet vectors are defined
unambiguously, while signs for edge vectors are defined ``locally" at each
vertex as described in Lemma~\ref{lm}. Our pictures contain a polytope (a
polygon in our case), facet and edge vectors, and an orientation of the
polytope (which is determined in our case by a cyclic order of vertices).

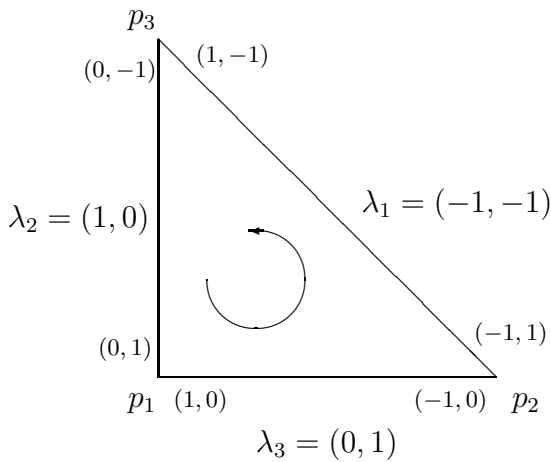
\begin{figure}
\begin{center}
\begin{picture}(100,60)
  \put(20,10){\line(0,1){45}}
  \put(20,10){\line(1,0){45}}
  \put(20,55){\line(1,-1){45}}
  \put(33,23){\oval(13,13)[b]}
  \put(33,23){\oval(13,13)[tr]}
  \put(34,29.5){\vector(-1,0){2}}
  \put(16,6){{\Large $p_1$}}
  \put(22,6){\small $(1,0)$}
  \put(54,6){\small $(-1,0)$}
  \put(67,6){{\Large $p_2$}}
  \put(33,0){{\Large $\lambda_3=(0,1)$}}
  \put(12,13){\small $(0,1)$}
  \put(0,30){{\Large $\lambda_2=(1,0)$}}
  \put(10,50){\small $(0,-1)$}
  \put(16,57){{\Large $p_3$}}
  \put(25,52){\small $(1,-1)$}
  \put(62,15){\small $(-1,1)$}
  \put(47,32){{\Large $\lambda_1=(-1,-1)$}}
\end{picture}%
\caption{$\tau(\C P^2)\oplus\C\simeq\bar\eta\oplus\bar\eta\oplus\bar\eta$}
\label{fig1}
\end{center}
\end{figure}

Figure 1 corresponds to $\C P^2$ regarded as a toric variety. Hence, a
stably complex structure is determined by the standard complex structure
on $\C P^2$, i.e. via the isomorphism of bundles $\tau(\C
P^2)\oplus\C\simeq\bar\eta\oplus\bar\eta\oplus\bar\eta$, where $\C$ is the
trivial complex line bundle and $\eta$ is the Hopf line bundle over $\C
P^2$. The orientation is defined by the complex structure. As we have
pointed out above, the toric variety $\C P^2$ arises from an integral
polytope (2-dimensional simplex with vertices $(0,0)$, $(1,0)$ and $(0,1)$
in this case). Facet vectors here are primitive normal vectors to facets
pointing inside the polytope, and edge vectors are primitive vectors
along edges pointing out of a vertex. This can be seen on Figure~1.
Let us calculate the Todd genus and the signature by
means of Corollary~\ref{signat} and Theorem~\ref{todd}. We have
$\sigma(p_1)=\sigma(p_2)=\sigma(p_3)=1$. Take $\nu=(1,2)$,
then $\ind(p_1)=0$, $\ind(p_2)=1$, $\ind(p_3)=2$
(remember that the index is the number of negative scalar products of edge
vectors with $\nu$). Thus, $\sign(\C P^2)= \sign(\C
P^2,\bar\eta\oplus\bar\eta\oplus\bar\eta)=1$, $\td(\C P^2)= \td(\C
P^2,\bar\eta\oplus\bar\eta\oplus\bar\eta)=1$.

\begin{figure}
\begin{center}
\begin{picture}(100,60)
  \put(20,10){\line(0,1){45}}
  \put(20,10){\line(1,0){45}}
  \put(20,55){\line(1,-1){45}}
  \put(33,23){\oval(13,13)[b]}
  \put(33,23){\oval(13,13)[tr]}
  \put(34,29.5){\vector(-1,0){2}}
  \put(16,6){{\Large $p_1$}}
  \put(22,6){\small $(-1,0)$}
  \put(54,6){\small $(1,0)$}
  \put(67,6){{\Large $p_2$}}
  \put(33,0){{\Large $\lambda_3=(0,1)$}}
  \put(12,13){\small $(0,1)$}
  \put(0,30){{\Large $\lambda_2{=}({-}1,0)$}}
  \put(10,50){\small $(0,-1)$}
  \put(16,57){{\Large $p_3$}}
  \put(25,52){\small $(-1,-1)$}
  \put(62,15){\small $(1,1)$}
  \put(47,32){{\Large $\lambda_1=(1,-1)$}}
\end{picture}%
\caption{$\tau(\bcp)\oplus\C\simeq\bar\eta\oplus\bar\eta\oplus\bar\eta$}
\end{center}
\end{figure}

Changing the orientation of polytope on Figure~1 causes changing of
signs of the vertices: $\sigma(p_1)=\sigma(p_2)=\sigma(p_3)=-1$, while the
indices of vertices remain unchanged. This corresponds to reversing the
canonical orientation of $\C P^2$. Denoting the resulting manifold $\bcp$, we
obtain $\sign(\bcp)= \sign(\bcp,\bar\eta\oplus\bar\eta\oplus\bar\eta)=-1$,
$\td(\bcp)= \td(\bcp,\bar\eta\oplus\bar\eta\oplus\bar\eta)=-1$. The same
stably complex structure can be obtained from the initial orientation of the
polytope by taking another facet vectors, as it is shown on Figure~2. In this
case one has $$ \sigma(p_1)=\begin{vmatrix} -1&0\\0&1 \end{vmatrix}=-1,\quad
  \sigma(p_2)=\begin{vmatrix} 1&1\\1&0 \end{vmatrix}=-1,\quad
  \sigma(p_3)=\begin{vmatrix} 0&-1\\-1&-1 \end{vmatrix}=-1.
$$
Again, we can take $\nu=(1,2)$, then $\ind(p_1)=1$, $\ind(p_2)=0$,
$\ind(p_3)=2$. Thus, we see again that $\sign(\bcp)=
\sign(\bcp,\bar\eta\oplus\bar\eta\oplus\bar\eta)=-1$, $\td(\bcp)=
\td(\bcp,\bar\eta\oplus\bar\eta\oplus\bar\eta)=-1$.

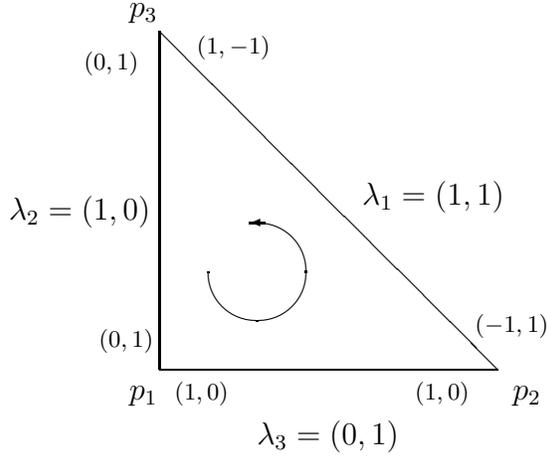
\begin{figure}
\begin{center}
\begin{picture}(100,60)
  \put(20,10){\line(0,1){45}}
  \put(20,10){\line(1,0){45}}
  \put(20,55){\line(1,-1){45}}
  \put(33,23){\oval(13,13)[b]}
  \put(33,23){\oval(13,13)[tr]}
  \put(34,29.5){\vector(-1,0){2}}
  \put(16,6){{\Large $p_1$}}
  \put(22,6){\small $(1,0)$}
  \put(54,6){\small $(1,0)$}
  \put(67,6){{\Large $p_2$}}
  \put(33,0){{\Large $\lambda_3=(0,1)$}}
  \put(12,13){\small $(0,1)$}
  \put(0,30){{\Large $\lambda_2=(1,0)$}}
  \put(10,50){\small $(0,1)$}
  \put(16,57){{\Large $p_3$}}
  \put(25,52){\small $(1,-1)$}
  \put(62,15){\small $(-1,1)$}
  \put(47,32){{\Large $\lambda_1=(1,1)$}}
\end{picture}%
\caption{$\tau(\C P^2)\oplus\C\simeq\eta\oplus\bar\eta\oplus\bar\eta$}
\label{fig3}
\end{center}
\end{figure}

Our third example (Figure~3) is $\C P^2$ with the standard orientation and
stably complex structure determined by the isomorphism
$\tau(\C P^2)\oplus\C\simeq\eta\oplus\bar\eta\oplus\bar\eta$
(this is obtained from Figure~1 by changing the sign of
first facet vector). Then we have
$$
  \sigma(p_1)=\begin{vmatrix} 1&0\\0&1 \end{vmatrix}=1,\quad
  \sigma(p_2)=\begin{vmatrix} -1&1\\1&0 \end{vmatrix}=-1,\quad
  \sigma(p_3)=\begin{vmatrix} 0&1\\1&-1 \end{vmatrix}=-1.
$$
Taking $\nu=(1,2)$, we find
$\ind(p_1)=0$, $\ind(p_2)=0$, $\ind(p_3)=1$. Thus,
$\sign(\C P^2,\eta\oplus\bar\eta\oplus\bar\eta)=1$,
$\td(\C P^2,\eta\oplus\bar\eta\oplus\bar\eta)=0$. Note that in this case
formula~(\ref{cn}) gives $c_n(\C
P^2,\eta\oplus\bar\eta\oplus\bar\eta)[\C P^2]=
\sigma(p_1)+\sigma(p_2)+\sigma(p_3)=-1$
(this could be also checked directly), while the Euler number of
$\C P^2$ is $c_n(\C P^2,\bar\eta\oplus\bar\eta\oplus\bar\eta)[\C P^2]=3$.

The author is grateful to Prof.~V.\,M.~Buchstaber for stimulating
discussions and useful recommendations.

\end{document}